\documentclass[12pt,twoside]{amsart}

\usepackage{amsmath,amsthm,amscd,amssymb,mathrsfs,graphicx,amsfonts,mathrsfs}
\usepackage{amsfonts}
\usepackage{amssymb,enumerate}
\usepackage{amsthm}
\usepackage[all]{xy}
\usepackage{hyperref}
\usepackage{needspace}
\hypersetup{pdftitle={Gorenstein and Ding injectivity over coherent rings},pdfauthor={Guocheng Dai and Jiangsheng Hu}}
\allowdisplaybreaks[4] \footskip=12pt
\renewcommand{\uppercasenonmath}[1]{}

\makeatletter
\renewcommand{\section}{\@startsection{section}{1}%
  \z@{.7\linespacing\@plus\linespacing}{.5\linespacing}%
  {\normalfont\bfseries\centering\let\@secnumfont\bfseries}}
\makeatother
\numberwithin{equation}{section} \theoremstyle{plain}
\newtheorem*{thm*}{Main Theorem}
\newtheorem{thm}{Theorem}[section]
\newtheorem{cor}[thm]{Corollary}
\newtheorem*{cor*}{Corollary}
\newtheorem{lem}[thm]{Lemma}
\newtheorem*{lem*}{Lemma}

\newtheorem*{fact*}{Fact}

\newtheorem*{nota*}{Notation}
\newtheorem{prop}[thm]{Proposition}
\newtheorem*{prop*}{Proposition}
\theoremstyle{remark}
\newtheorem{rem}[thm]{Remark}
\newtheorem*{rem*}{Remark}
\theoremstyle{plain}

\newtheorem*{observation*}{Observation}

\newtheorem*{exa*}{Example}

\newtheorem*{df*}{Definition}

\newtheorem*{con*}{Construction}

\newcommand{\Hom}{\operatorname{Hom}}

\newcommand{\ann}{\operatorname{Ann}}

\newcommand{\im}{\operatorname{Im}}

\renewcommand{\geq}{\geqslant}
\renewcommand{\leq}{\leqslant}
\renewcommand{\ker}{\Ker}

\newcommand{\GI}{\mathcal{GI}}
\newcommand{\DI}{\mathcal{DI}}

\newcommand{\Inj}{\operatorname{Inj}}
\newcommand{\Ext}{\operatorname{Ext}}
\renewcommand{\im}{\operatorname{im}}
\newcommand{\Cof}{\operatorname{Cof}}
\newcommand{\Gid}{\operatorname{Gid}}
\newcommand{\Did}{\operatorname{Did}}
\newcommand{\Spec}{\operatorname{Spec}}
\renewcommand{\ann}{\operatorname{ann}}
\newcommand{\FPid}{\operatorname{FP\!-\!id}}
\newcommand{\Loc}{\operatorname{Loc}}
\renewcommand{\ker}{\operatorname{ker}}
\begin{document}
\begin{center}
{\large  \bf Gorenstein and Ding injectivity over coherent rings}

\vspace{0.5cm} Guocheng Dai and Jiangsheng Hu\\
\end{center}

\bigskip
\centerline { \bf  Abstract}
%\bigskip
\leftskip10truemm \rightskip10truemm\noindent
For a commutative Ding--Chen ring $B$, we prove that Gorenstein
and Ding injectivity coincide over $R=B[x,y]/(x,y)^2$ if and only if
injective $B$-modules generate the unbounded derived category
of $B$. We construct a commutative von Neumann regular
ring $B$ whose unbounded derived category is not generated by
its pure injective objects. For this $B$, the coherent ring $R$
admits a strongly Gorenstein injective module $G$ which is not
Ding injective. For every prime ideal $\mathfrak p$ of $R$, the
$R_{\mathfrak p}$-module $G_{\mathfrak p}$ has infinite Gorenstein
injective dimension.
The ring $R$ can be chosen to have cardinality
$\aleph_\omega$, the least possible among coherent rings over
which Gorenstein and Ding injectivity differ.

\leftskip10truemm \rightskip10truemm \noindent
\vbox to 0.3cm{}

\leftskip0truemm \rightskip0truemm

\section{Introduction}

Gorenstein injectivity is defined using exact complexes of
injective modules which remain exact under Hom from injective
modules (see \cite{EJ95}). Ding injectivity requires the same
condition for all absolutely pure modules (see
\cite[p.~491]{MD08} and \cite[Definition 3.2]{Gil10}).
Write $\Inj(R)$, $\GI(R)$ and $\DI(R)$ for the classes of
injective, Gorenstein injective and Ding injective left
$R$-modules, respectively. Then
$\Inj(R)\subseteq\DI(R)\subseteq\GI(R)$, with
$\DI(R)=\GI(R)$ over left Noetherian rings.
The definitions and basic facts are recalled in
Section~\ref{sec:comparison}.

Gillespie asked whether $\DI(R)=\GI(R)$ over Ding--Chen
rings (see \cite[Section 3]{Gil10}); the two classes coincide
over these rings (see \cite[Corollary 4.3]{Gil17}).
Wang, Liu and Yang construct a Gorenstein injective module
which is not Gorenstein AC-injective over a noncoherent ring
(see \cite[Section 2.3 and Examples 24--25]{WLY21}).
Gorenstein AC-injectivity agrees with Ding injectivity over
coherent rings, but their example lies outside this setting.
The question considered here is whether every Gorenstein
injective module over a coherent ring must be Ding injective.

We study this question for the square-zero extensions
$R=B[x,y]/(x,y)^2$, where $B$ is a commutative Ding--Chen
ring. For a field $k$, Zhang constructs an acyclic one-periodic
free complex over $k[x,y]/(x,y)^2$
(see \cite[Lemma 2.5]{ZhaPeriodic}).
We use the injective $R$-modules $\Hom_B(R,E)$, with $E$
injective over $B$, to adapt this construction.
The construction starts from a nonexact complex $U$ of
injective $B$-modules such that $\Hom_B(E,U)$ is exact for
every injective $B$-module $E$. Such a complex cannot exist
if injective $B$-modules generate the unbounded derived
category $\mathbf D(B)$. This leads to the following
characterization of $\DI(R)=\GI(R)$. For a class $\mathcal S$ of objects of
$\mathbf D(B)$, let $\Loc(\mathcal S)$ denote the localizing
subcategory it generates.

\Needspace{10\baselineskip}
\begin{thm}\label{thm:regular-comparison}
Let $B$ be a commutative Ding--Chen ring, and set
$R=B[x,y]/(x,y)^2$. The following conditions are equivalent:
\begin{enumerate}
\renewcommand{\theenumi}{\roman{enumi}}
\renewcommand{\labelenumi}{\textup{(\theenumi)}}
\setlength{\itemsep}{2pt}
\setlength{\parsep}{0pt}
\item\label{cond:equality}
$\DI(R)=\GI(R)$.
\item\label{cond:exactness}
Every complex $U$ of injective $B$-modules is exact whenever
$\Hom_B(E,U)$ is exact for every injective $B$-module $E$.
\item\label{cond:generation}
$\Loc(\Inj(B))=\mathbf D(B)$.
\end{enumerate}
\end{thm}

Condition~\textup{(\ref{cond:exactness})} says that the functors
$\Hom_B(E,-)$, with $E$ injective, detect exactness of complexes
of injective $B$-modules. Condition~\textup{(\ref{cond:generation})}
places this property in the study of generation by injective
modules (see \cite[Section 3]{Ric19}).
We obtain two applications by constructing a commutative
von Neumann regular algebra for which
condition~\textup{(\ref{cond:exactness})} fails.

The first application concerns generation by pure injective
objects in the triangulated sense
(see \cite[Section 2.3]{BHM26}). Breaz, Hrbek and Modoi discuss
whether these objects generate the derived category of every
commutative ring. Over a von Neumann regular ring, this is
equivalent to generation by injective modules
(see \cite[Remark 6.12]{BHM26}).

\begin{cor}\label{cor:pure-injective-generation}
For every field $k$, there is a commutative von Neumann
regular $k$-algebra $B$ of $k$-dimension $\aleph_\omega$
such that $\Loc(\mathcal P)\ne\mathbf D(B)$, where
$\mathcal P$ is the class of pure injective objects of
$\mathbf D(B)$.
\end{cor}

Corollary~\ref{cor:pure-injective-generation} answers negatively
the generation question discussed above.
For countable $k$, the ring $B$ has the least possible
cardinality $\aleph_\omega$ among commutative rings whose
derived categories are not generated by pure injective objects.
The lower bound follows from \cite[Remark 6.11(iii)]{BHM26}.

To state the second application, write $\Lambda_k:=k[x,y]/(x,y)^2$.
For strongly Gorenstein injective modules, the defining complex
can be chosen one-periodic (see \cite{BM07}).

\begin{cor}\label{cor:separation-example}
For every field $k$, there exist a commutative coherent
$k$-algebra $R$ of $k$-dimension $\aleph_\omega$
and a strongly Gorenstein injective $R$-module $G$
which is not Ding injective.
Moreover, one has $\DI(R)=\Inj(R)$.
For every prime ideal $\mathfrak p$ of $R$, one has 
$R_{\mathfrak p}\cong\Lambda_k$ and
$\Gid_{R_{\mathfrak p}}G_{\mathfrak p}=\infty$.
\end{cor}

Corollary~\ref{cor:separation-example} shows that coherence alone
does not suffice for Gorenstein and Ding injectivity to coincide.
All local rings in this example are Noetherian, so the two
classes coincide over each of them, although $G_{\mathfrak p}$
belongs to neither class. Couchot's example also shows that
Gorenstein injectivity can fail under localization at a
multiplicative set (see \cite[Example 1]{Cou06}); in our example,
for every prime ideal $\mathfrak p$ of $R$, the
$R_{\mathfrak p}$-module $G_{\mathfrak p}$ has infinite Gorenstein
injective dimension.

For countable $k$, the ring in
Corollary~\ref{cor:separation-example} has the least possible
cardinality $\aleph_\omega$ among coherent rings over which
Gorenstein and Ding injectivity differ. Indeed, every absolutely pure
module over a left coherent ring of cardinality at most
$\aleph_n$ has injective dimension at most $n+1$ for $n<\omega$
(see \cite[Proposition 2.3]{Pos17}). Dimension shifting
therefore proves $\GI(R)=\DI(R)$ for every left coherent
ring $R$ of cardinality less than $\aleph_\omega$.

The two applications use the same construction. We construct
a commutative von Neumann regular $k$-algebra $B$ such that
$\Ext_B^m(E,B)=0$ for every injective $B$-module $E$ and
every $m\ge0$. An injective resolution $U$ of $B$ in
nonnegative degrees has $H^0(U)=B\ne0$, whereas
$\Hom_B(E,U)$ is exact for every injective $B$-module $E$.
Thus $U$ violates condition~\textup{(\ref{cond:exactness})} of
Theorem~\ref{thm:regular-comparison}.

%The construction sends $1\in B$ to an Ext class that
%remains nonzero after localization at every prime.
%Periodicity and dimension shifting establish the
%infinite local Gorenstein injective dimensions.

The paper is organized as follows.
Section~\ref{sec:comparison} recalls the main definitions and
basic facts used in the proofs. Section~\ref{sec:transfer}
proves Theorem~\ref{thm:regular-comparison}.
Section~\ref{sec:applications} contains two applications:
Subsection~\ref{sec:generation} constructs the ring used in
both applications and proves
Corollary~\ref{cor:pure-injective-generation};
Subsection~\ref{sec:realization} proves
Corollary~\ref{cor:separation-example} and derives an example
with both strict inclusions $\Inj\subsetneq\DI\subsetneq\GI$.

\section{Preliminaries}
\label{sec:comparison}

All rings are associative with identity and nonzero. Modules
are left modules unless otherwise stated. Complexes are cochain
and may be unbounded; a module is regarded as a complex in degree
zero. We write $Z^nJ$ and $H^nJ$ for the cycles and cohomology of
a complex $J$.
For complexes $P$ and $U$, we write $\Hom_R(P,U)$ for the
Hom complex of abelian groups, with
$\Hom_R(P,U)^n=\prod_{a\in\mathbb Z}\Hom_R(P^a,U^{a+n})$
and differential $d(f)=d_Uf-(-1)^nf d_P$ for $f$ of degree $n$.
We write $\Gid$, $\Did$ and $\FPid$ for Gorenstein injective,
Ding injective and FP-injective dimensions, respectively.

\subsection{Gorenstein and Ding injective modules}
\label{subsec:injective-modules}
A left $R$-module $N$ is \emph{absolutely pure}, or
\emph{FP-injective}, if $\Ext_R^1(F,N)=0$ for every finitely
presented left $R$-module $F$ (see \cite[Definition 3.1]{Gil10}).
An exact complex $I$ of injective $R$-modules is
\emph{totally acyclic} if $\Hom_R(E,I)$ is exact for every
injective $R$-module $E$, and \emph{Ding acyclic} if
$\Hom_R(A,I)$ is exact for every absolutely pure $R$-module $A$.
A module $M$ is \emph{Gorenstein injective} if $M=Z^0I$
for a totally acyclic complex $I$; such a complex is a
\emph{complete injective resolution} of $M$ (see \cite{EJ95}).
A cycle of a Ding acyclic complex is called \emph{Ding
injective} (see \cite[p.~491]{MD08} and
\cite[Definition 3.2]{Gil10}).
A complex is \emph{one-periodic} if it has the same module
and differential in every degree. A module is \emph{strongly
Gorenstein injective} if it has a one-periodic complete
injective resolution (see \cite{BM07}).
The inclusions $\Inj(R)\subseteq\DI(R)\subseteq\GI(R)$
follow from these definitions. Over a left Noetherian ring,
absolutely pure modules are injective, so $\DI(R)=\GI(R)$.

For a Ding acyclic complex $J$, the complex
$\Hom_R(L,J)$ is exact for every $R$-module $L$
of finite FP-injective dimension.
Indeed, if $\FPid_R L=d$, then the $d$-th cosyzygy
of $L$ in an injective resolution is absolutely pure,
by dimension shifting against finitely presented modules.
Induction along the resolution proves the assertion,
since Hom into each term of $J$ is exact on short
exact sequences.

A ring is called a \emph{Ding--Chen ring} if it is two-sided
coherent and has finite self FP-injective dimension on both
sides (see \cite[Definition 4.1]{Gil10}).
By \cite[Corollary 4.3]{Gil17}, every exact complex of injective
modules over a Ding--Chen ring is Ding acyclic.

For a class $\mathcal C$ of $R$-modules, we use
${}^{\perp}\mathcal C$ for the modules $L$ such that
$\Ext_R^1(L,C)=0$ for every $C\in\mathcal C$.

\subsection{Derived categories}
\label{subsec:complexes}

We write $\mathbf D(R)$ for the unbounded derived category of
left $R$-modules. A complex $I$ is \emph{K-injective} if
$\Hom_R(A,I)$ is exact for every exact complex $A$
(see \cite{Spa88}).
For a class $\mathcal S$ of objects of $\mathbf D(R)$,
$\Loc(\mathcal S)$ denotes the smallest full triangulated
subcategory containing $\mathcal S$ and closed under arbitrary
coproducts. Its right orthogonal consists of the objects $X$
such that $\Hom_{\mathbf D(R)}(S,X[n])=0$ for all
$S\in\mathcal S$ and $n\in\mathbb Z$.
For an abelian group $M$, write
$M^+=\Hom_{\mathbb Z}(M,\mathbb Q/\mathbb Z)$.

A morphism $f:X\to Y$ in $\mathbf D(R)$ is a \emph{pure
monomorphism} if $\Hom_{\mathbf D(R)}(C,f)$ is injective for
every compact object $C$.
An object $P$ is \emph{pure injective} if every pure
monomorphism $P\to Y$ splits. We use this triangulated
notion of purity (see \cite[Section 2.3]{BHM26}).

\subsection{Square-zero extensions}
\label{subsec:square-zero}

Write $\Lambda_k=k[x,y]/(x,y)^2$ for a field $k$.
For a commutative ring $B$, set $R=B[x,y]/(x,y)^2$ and write
$\mathfrak r=(x,y),\; \Cof_B(L)=\Hom_B(R,L),\; F_B=\Cof_B(B)$.
For a ring map $B\to A$ and an $A$-module $N$, the
notation $N_B$ denotes restriction of scalars to $B$;
the same convention applies degreewise to complexes.
A $B$-module is inflated to $R$ through the quotient
map $R\to B=R/\mathfrak r$.
The $R$-action on $\Cof_B(L)$ is defined by
$(rf)(s)=f(sr)$.
The functor $\Cof_B$ is right adjoint to restriction of scalars
and preserves injectives. It is exact since $R$ is free over $B$.
If $B$ is an algebra over a field $k$, then
$R=B\otimes_k\Lambda_k$ and
$F_B\cong B\otimes_k\Hom_k(\Lambda_k,k)$.
Here $\Lambda_k$ is the path algebra of the quiver with
one vertex and two loops, modulo all paths of length two.

Using the $B$-basis $1,x,y$ of $R$, we identify
$\Cof_B(L)$ with $L^3$.
For the dual basis $\lambda_0,\lambda_x,\lambda_y$,
the actions of $x$ and $y$ are
$x(a,b,c)=(b,0,0),\; y(a,b,c)=(c,0,0)$.
In particular, as $B$-modules,
$R/\mathfrak r\cong B$ and $\mathfrak r\cong B^2$,
whereas $F_B/\mathfrak rF_B\cong B^2$ and
$\mathfrak rF_B\cong B$.

\begin{samepage}
\begin{lem}\label{lem:finite-detector}
Let $B$ be a commutative ring, set $R=B[x,y]/(x,y)^2$,
and write $F_B=\Hom_B(R,B)$. The $R$-module $F_B$ is finitely presented and
$\FPid_R F_B=\FPid_B B$. If $B$ is coherent, so is $R$.
If $B$ is von Neumann regular, then $F_B$ is absolutely
pure. In this case, for $\mathfrak p\in\Spec R$,
$\mathfrak q=\mathfrak p\cap B$ and $k(\mathfrak q)=B/\mathfrak q$, one has 
$$
 R_{\mathfrak p}\cong\Lambda_{k(\mathfrak q)},\;
 (F_B)_{\mathfrak p}\cong
 \Hom_{k(\mathfrak q)}(\Lambda_{k(\mathfrak q)},k(\mathfrak q)).
$$
In particular, $(F_B)_{\mathfrak p}$ is injective over $R_{\mathfrak p}$.
\end{lem}
\end{samepage}
\begin{proof}
Suppose $B$ is coherent. The ring $R$ is finite free
over $B$.
The kernel of a map between finitely generated free
$R$-modules is therefore finitely generated over $B$,
and hence over $R$.
This proves coherence. For arbitrary $B$, the coordinate
actions give the finite presentation
$$
 R^3\xrightarrow{\left(\begin{smallmatrix}y&0&x\\0&x&-y\end{smallmatrix}\right)}
 R^2\longrightarrow F_B\longrightarrow0.
$$
The map $R^2\to F_B$ sends the standard basis to
$\lambda_x,\lambda_y$.
Indeed, its kernel is generated by
$(y,0),(0,x),(x,-y)$.
The functor $\Cof_B$ is exact and preserves injectives.
By adjunction, for a left $R$-module $M$ we have
$\Ext_R^q(M,F_B)\cong\Ext_B^q(M_B,B)$.
Restriction takes finitely presented $R$-modules to
finitely presented $B$-modules. Conversely, for a
finitely presented $B$-module $L$, the induced module
$R\otimes_B L$ is finitely presented over $R$ and
restricts to $L^3$. This proves the equality of
FP-injective dimensions.

Assume now that $B$ is von Neumann regular. Its finitely
presented modules are projective, so $F_B$ is absolutely pure.
Every prime ideal of $R$ contains $\mathfrak r$. Since $B$ is commutative von Neumann regular,
$B_{\mathfrak q}\cong B/\mathfrak q$.
The two displayed identifications follow by localizing
the three basis elements and the actions of $x$ and $y$.
The vector-space dual of a finite-dimensional algebra over
a field is injective by the same adjunction.
\end{proof}

\section{Proof of Theorem~\ref{thm:regular-comparison}}
\label{sec:transfer}

Throughout this section, let $B$ be a commutative ring and
set $R=B[x,y]/(x,y)^2$. 
We first compare the cohomology of a complex $J$ of injective
$R$-modules with that of $\Hom_R(B,J)$.
We then adapt \cite[Lemma 2.5]{ZhaPeriodic} to construct a
periodic complex of injective $R$-modules and complete the
proof of Theorem~\ref{thm:regular-comparison}.

We next describe injective $R$-modules and two associated
short exact sequences. For an $R$-module $N$, write
$S(N)=\{z\in N:xz=yz=0\}$. Evaluation at $1$ identifies
$\Hom_R(B,N)$ with $S(N)$. For complexes, $S$ is applied degreewise.

\begin{lem}\label{lem:injective-normal-form}
Every injective $R$-module is of the form $\Cof_B(E)$
for an injective $B$-module $E$. If $J$ is a complex
of injective $R$-modules and $U=S(J)$, then $U$ has
injective $B$-terms and there is a short exact sequence
of complexes of $B$-modules, split in each degree,
$$
 0\longrightarrow U\longrightarrow J_B\longrightarrow U\oplus U
   \longrightarrow0.
$$
For every $B$-module $L$, set
$C_L=\Hom_B(L,U)$ and
$Y_L=\Hom_R(\Cof_B(L),J)$. There is also a short exact sequence of complexes of
$B$-modules, split in each degree,
$$
 0\longrightarrow C_L\oplus C_L\longrightarrow Y_L\longrightarrow C_L
   \longrightarrow0.
$$
The connecting maps of these sequences, after applying
$\Hom_B(L,-)$ to the first, are respectively
$$
 (v_x,v_y):H^n(C_L)^2\longrightarrow H^{n+1}(C_L),
 \;
 \binom{v_x}{v_y}:H^n(C_L)\longrightarrow H^{n+1}(C_L)^2,
$$
for the same two maps $v_x,v_y$. If $U$ is contractible,
then $J$ is contractible over $R$.
\end{lem}
\begin{proof}
For an injective $R$-module $Q$, the module $S(Q)$
is injective over $B$, by adjunction with restriction of
scalars along the quotient map $R\to B$.
The submodule $S(Q)$ is essential:
any element not annihilated by $\mathfrak r$ generates a
nonzero element annihilated by $\mathfrak r$, since
$\mathfrak r^2=0$. The same is true of the copy of
$S(Q)$ inside the injective module $\Cof_B(S(Q))$.
By uniqueness of injective envelopes, there is an isomorphism
$Q\cong\Cof_B(S(Q))$ fixing this submodule.

Choose these identifications in each degree. The coordinate
actions show that each $R$-linear differential of $J$
has the form
$$
 d_J^n=\begin{pmatrix}u^n&v_x^n&v_y^n\\0&u^n&0\\0&0&u^n\end{pmatrix},
 \; u^{n+1}v_i^n+v_i^{n+1}u^n=0.
$$
Here $u$ is the differential of $U$. Inclusion of the
first coordinate and projection onto the last two give
the first short exact sequence. An $R$-linear map
$\Cof_B(L)\to\Cof_B(U^n)$ has the same upper triangular
form, determined by three maps $\alpha,\beta_x,\beta_y$.
The Hom differential is
$$
 (\alpha,\beta_x,\beta_y)\longmapsto
 (u\alpha,u\beta_x+v_x\alpha,u\beta_y+v_y\alpha).
$$
Including the two $\beta$-coordinates and projecting onto
$\alpha$ gives the second short exact sequence.
The connecting maps of the two sequences are the stated
row and column, respectively.

Finally, if $h$ contracts $U$, set
$a_i^n=h^{n+1}v_i^n$. Then
$v_i^n=u^na_i^n-a_i^{n+1}u^n$.
The change of coordinates
$(a,b,c)\mapsto(a+a_x^nb+a_y^nc,b,c)$ identifies $J$
with $\Cof_B(U)$, which is contracted by $\Cof_B(h)$.
\end{proof}

\begin{prop}\label{prop:single-source}
Let $J$ be an exact complex of injective $R$-modules,
and write $U=S(J)$. For every $n\in\mathbb Z$, there is
an isomorphism of $R$-modules
$H^n\Hom_R(F_B,J)\cong H^{n-1}(U)^3$,
where the action on the right factors through $R\to B$.
If $B$ is von Neumann regular, then $J$ is Ding acyclic
if and only if $U$ is exact. These conditions are
equivalent to contractibility of $J$, and
$\DI(R)=\Inj(R)$.
\end{prop}
\begin{proof}
With $L=B$ in Lemma~\ref{lem:injective-normal-form}, write
$V^n=H^n(U)$. Exactness of $J_B$ makes each row map
$(V^{n-1})^2\to V^n$ an isomorphism, and hence each
column map is injective. The cohomology sequence of the
second short exact sequence in
Lemma~\ref{lem:injective-normal-form} identifies
$H^n\Hom_R(F_B,J)$ with the cokernel of
$V^{n-1}\to(V^n)^2$. In the coordinates given by the
inverse of the row isomorphism, this column sends
$z$ to $(z,0,0,z)$ in $(V^{n-1})^4$. Its cokernel
is identified with $(V^{n-1})^3$ by
$(a,b,c,d)\mapsto(b,c,d-a)$. The maps from the two
$\beta$-coordinates have
image in $S(J)$, on which $x$ and $y$ act by zero.
Thus this isomorphism is $R$-linear when
$R$ acts on its right-hand side through $R\to B$.
Since the row maps are isomorphisms in every degree,
vanishing of $H^n\Hom_R(F_B,J)$ for one
$n$ forces $U$ to be exact.

Assume now that $B$ is von Neumann regular.
The module $F_B$ is absolutely pure by
Lemma~\ref{lem:finite-detector}. Thus, if $J$ is Ding
acyclic, $U$ is exact.
Conversely, if $U$ is exact, it is contractible by
\cite[Theorem 4.1(1)]{BCE}, since every $B$-module is absolutely
pure. Thus $J$ is contractible
by Lemma~\ref{lem:injective-normal-form}. A contractible
complex is Ding acyclic, so these conditions are equivalent.
In particular, every Ding injective module is injective.
The reverse inclusion follows from the definition, and hence
$\DI(R)=\Inj(R)$.
\end{proof}

\begin{lem}\label{lem:ding-chen-criterion}
Let $B$ be a commutative Ding--Chen ring, let
$R=B[x,y]/(x,y)^2$, and let $J$ be an exact complex
of injective $R$-modules. Then $J$ is Ding acyclic
if and only if $S(J)$ is exact.
\end{lem}
\begin{proof}
Write $U=S(J)$.
If $J$ is Ding acyclic, then
$\Hom_R(F_B,J)$ is exact: $F_B$ has finite
FP-injective dimension over $R$ by Lemma~\ref{lem:finite-detector},
so the dimension-shifting argument in Subsection~\ref{subsec:injective-modules} applies.
By Proposition~\ref{prop:single-source}, $U$ is exact.

Conversely, suppose $U$ is exact. Let $A$ be absolutely pure
over $R$ and set $L=S(A)$. The ring $R=B\ltimes B^2$ is
coherent by Lemma~\ref{lem:finite-detector}, and $B^2$ is finitely
presented over $B$. By \cite[Lemma 3.1(1)]{Mao23}, $L$ is
absolutely pure over $B$ and there is an exact sequence
$$
 0\longrightarrow L\longrightarrow A
 \xrightarrow{a\mapsto(xa,ya)}L^2\longrightarrow0,
$$
where the end terms are inflated from $B$. The complex
$U$ is exact with injective $B$-terms. Since $L$ is
absolutely pure, $\Hom_B(L,U)$ is exact by
\cite[Corollary 4.3]{Gil17}. Applying
$\Hom_R(-,J)$ to the displayed sequence
proves exactness of $\Hom_R(A,J)$: the
complexes on the end terms identify with
$\Hom_B(L,U)$ and its direct sum with itself.
\end{proof}

\begin{rem}\label{rem:ding-injective-test}
When $B$ is Ding--Chen, Proposition~\ref{prop:single-source}
and Lemma~\ref{lem:ding-chen-criterion} show that a cycle $G$
of an exact complex of injective $R$-modules is Ding injective
if and only if $\Ext_R^m(F_B,G)=0$ for some $m\geq1$.
\end{rem}

The next construction proves
\textup{(\ref{cond:equality})}$\Rightarrow$\textup{(\ref{cond:exactness})}
in Theorem~\ref{thm:regular-comparison} by contraposition and
provides the Ext class used in Subsection~\ref{sec:realization}.

\begin{prop}
\label{prop:periodic-transfer}
Let $B$ be a commutative ring, set $R=B[x,y]/(x,y)^2$, and
write $F_B=\Hom_B(R,B)$. Suppose $U$ is a complex
of injective $B$-modules for which
$\Hom_B(E,U)$ is exact for every injective $B$-module $E$.
There is a strongly Gorenstein injective $R$-module $G$
and an exact sequence
$$
 0\longrightarrow G\longrightarrow I\longrightarrow G
   \longrightarrow0
$$
with $I$ injective and $\Ext_R^1(Q,G)=0$ for every injective
$R$-module $Q$.
For this $G$, there is a $B$-linear monomorphism
$H^0(U)\longrightarrow\Ext_R^1(F_B,G)$ whose image is
annihilated by $(x,y)$.
If $\FPid_B B<\infty$ and $H^0(U)\ne0$, then $G$ is not
Ding injective.
More precisely, set $W=\prod_{j\in\mathbb Z,\;r\ge0}H^j(U)$.
There are $R$-module isomorphisms
$\Ext_R^n(F_B,G)\cong W^3\;(n\ge1)$,
where $R$ acts on $W^3$ through $R\to B$.
Under the isomorphism for $n=1$, the monomorphism is the
inclusion of the $H^0(U)$-coordinate with $r=0$ in the
first copy of $W$.
\end{prop}

\begin{proof}
We first construct an exact one-periodic complex of injective
$R$-modules. The integer index records the degrees of $U$,
while the additional countable copies allow us to identify
$M^2$ with $M$. Define
$$
 M=\prod_{n\in\mathbb Z,\ r\ge0}U^n,\;
 (um)_{n,r}=d_U^{n-1}m_{n-1,r},\;
 (\tau m)_{n,r}=(-1)^n m_{n,r}.
$$
Define $s_x,s_y:M\to M$ by inserting the coordinates of
$m\in M$ into the even and odd positions of the second
index, respectively:
$$
 \begin{aligned}
 (s_xm)_{n,2r}&=m_{n,r},& (s_xm)_{n,2r+1}&=0,\\
 (s_ym)_{n,2r}&=0,& (s_ym)_{n,2r+1}&=m_{n,r}.
 \end{aligned}
$$
Set $v_x=\tau s_x$ and $v_y=\tau s_y$. The signs ensure
that $u$ anticommutes with $v_x$ and $v_y$. The row map
$(v_x,v_y):M^2\to M$ is an isomorphism, with inverse
$(p_x,p_y)$ obtained by applying $\tau$ and extracting
the even and odd coordinates. We have
$$
 \begin{gathered}
 u^2=0,\; uv_i+v_i u=0,\; p_i u=-up_i,\\
 p_i v_j=\delta_{ij},\; v_xp_x+v_yp_y=1.
 \end{gathered}
$$
The module $M$ is injective, being a product of injective
$B$-modules, so $I=\Cof_B(M)=M^3$ is injective over $R$.
Define
$$
 D=\begin{pmatrix}u&v_x&v_y\\0&u&0\\0&0&u\end{pmatrix}.
$$
The coordinate actions show that $D$ is $R$-linear,
and the identities above imply $D^2=0$. The $B$-linear map
$H(a,b,c)=(0,p_xa,p_ya)$ satisfies $DH+HD=1$.
Thus the complex $J$, with $J^n=I$ and $d_J^n=D$ in every
degree, is contractible as a complex of $B$-modules, and hence
exact as a complex of $R$-modules. Set $G=Z^0J$. By periodicity,
there is an exact sequence $0\to G\to I\to G\to0$, and the
nonnegative part of $J$ is an injective resolution of $G$.

Write $V=S(J)$, the one-periodic complex with term
$M$ and differential $u$. For an injective $B$-module
$E$, the identification
$\Hom_B(E,M)\cong \prod_{j\in\mathbb Z,\;r\ge0}\Hom_B(E,U^j)$
and exactness of products show that
$\Hom_B(E,V)$ is exact.
By the second short exact sequence of
Lemma~\ref{lem:injective-normal-form},
$\Hom_R(\Cof_B(E),J)$ is exact.
Every injective $R$-module has this form, by the same
lemma. Thus $J$ is totally acyclic, so $G$ is strongly
Gorenstein injective.

Exactness of products also identifies $H^q(V)$ with
$W$ in every degree $q$. Since the nonnegative part
of $J$ is an injective resolution of $G$,
the asserted Ext isomorphisms follow from
Proposition~\ref{prop:single-source}, including the stated
$R$-action. Inclusion of the $H^0(U)$-coordinate
with second index zero in the first copy of $W$ induces
the stated monomorphism.
If $\FPid_B B<\infty$, then $F_B$ has finite
FP-injective dimension over $R$ by Lemma~\ref{lem:finite-detector}.
Ding injectivity of $G$ would then force $\Ext_R^1(F_B,G)=0$,
contrary to the Ext isomorphism when $H^0(U)\ne0$.
\end{proof}

We are now in a position to prove Theorem~\ref{thm:regular-comparison} in the introduction.

\begin{proof}[\bf Proof of Theorem~\ref{thm:regular-comparison}]
\textup{(\ref{cond:equality})}$\Rightarrow$\textup{(\ref{cond:exactness})}.
Assume $\DI(R)=\GI(R)$, and let $U$ be a complex of injective
$B$-modules such that $\Hom_B(E,U)$ is exact for every
injective $B$-module $E$. Suppose that $U$ is not exact.
Choose $t\in\mathbb Z$ with $H^t(U)\ne0$. The complex $U[t]$
has injective terms, satisfies the same Hom exactness condition,
and has $H^0(U[t])=H^t(U)\ne0$. Since $B$ is Ding--Chen,
$\FPid_B B<\infty$. By Proposition~\ref{prop:periodic-transfer},
applied to $U[t]$, there is a strongly Gorenstein injective
$R$-module which is not Ding injective. This contradicts
\textup{(\ref{cond:equality})}. Hence $U$ is exact.

\textup{(\ref{cond:exactness})}$\Rightarrow$\textup{(\ref{cond:equality})}.
Let $G\in\GI(R)$, and choose a totally acyclic complex $J$
with $G=Z^0J$. Write $U=S(J)$. By
Lemma~\ref{lem:injective-normal-form}, $U$ has injective
$B$-terms. We show that $\Hom_B(E,U)$ is exact for every
injective $B$-module $E$.
For such a module $E$, the module $\Cof_B(E)$ is injective
over $R$, so $\Hom_R(\Cof_B(E),J)$ is exact.
Also $J_B$ is an exact complex of injective $B$-modules.
Since $B$ is Ding--Chen, \cite[Corollary 4.3]{Gil17}
shows that $\Hom_B(E,J_B)$ is exact.
The long exact cohomology sequences associated with the two
short exact sequences in Lemma~\ref{lem:injective-normal-form}
show that both connecting maps are isomorphisms when $L=E$.
Write $V_E^n=H^n\Hom_B(E,U)$.
Injectivity of the row implies that
$v_y:V_E^n\to V_E^{n+1}$ is injective.
Surjectivity of the column expresses every
$(w,0)\in(V_E^{n+1})^2$ as $(v_xz,v_yz)$.
Then $z=0$, so $w=0$. Thus $\Hom_B(E,U)$ is exact.
Since $E$ was arbitrary, \textup{(\ref{cond:exactness})}
implies that $U$ is exact. By
Lemma~\ref{lem:ding-chen-criterion}, $J$ is Ding acyclic,
and hence $G$ is Ding injective. As $G$ was arbitrary,
$\GI(R)\subseteq\DI(R)$. The reverse inclusion always
holds, so \textup{(\ref{cond:equality})} follows.

\textup{(\ref{cond:exactness})}$\Longleftrightarrow$\textup{(\ref{cond:generation})}.
Let $U$ be a complex of injective $B$-modules. By
\cite[Proposition 3.11]{Spa88}, choose a quasi-isomorphism
$U\to U'$, where $U'$ is K-injective with injective
terms. The cone $C$ of this quasi-isomorphism is exact with
injective terms. By \cite[Corollary 4.3]{Gil17}, $\Hom_B(E,C)$
is exact for every injective $B$-module $E$. Hence
$\Hom_B(E,U)\to\Hom_B(E,U')$ is a quasi-isomorphism, and
$\Hom_B(E,U)$ computes $\mathbf R\Hom_B(E,U)$.
Thus a complex $U$ of injective terms for which
$\Hom_B(E,U)$ is exact for every injective $B$-module $E$
represents an object in the right orthogonal to $\Inj(B)$.
Conversely, every object in that orthogonal admits a
K-injective representative $V$ with injective terms, and
$\Hom_B(E,V)$ is exact for every injective $B$-module $E$.
Consequently, \textup{(\ref{cond:exactness})} holds if and only
if this right orthogonal is zero.

To relate this vanishing to generation, we show that
$\Loc(\Inj(B))$ is generated by a single module. Set $Q=B^{++}$.
The module $B^+$ is injective and has finite flat dimension by
\cite[Theorem 4.2]{Gil10}; dualizing a finite flat resolution
shows that $Q$ has finite injective dimension. Since $B$ is
coherent, character duality also implies that $Q$ is flat
(see \cite[proof of Proposition 3.14]{Gil10}).
The evaluation map $B\to Q$ is pure, so each injective
$B$-module $E$ is a summand of $Q\otimes_B E$.
Tensoring with the flat module $Q$ preserves triangles
and coproducts; since $\Loc(B)=\mathbf D(B)$, it follows
that $Q\otimes_B E\in\Loc(Q)$.
Conversely, a finite injective resolution of $Q$ places
$Q$ in $\Loc(\Inj(B))$. Therefore
$\Loc(\Inj(B))=\Loc(Q)$.
By the Bousfield localization argument in
\cite[proof of Theorem 4.4]{Ric19}, each object of
$\mathbf D(B)$ fits into a triangle with first term in
$\Loc(Q)$ and last term in its right orthogonal.
That orthogonal is zero exactly when
$\Loc(Q)=\mathbf D(B)$. Since $\Loc(\Inj(B))=\Loc(Q)$,
this proves
\textup{(\ref{cond:exactness})}$\Longleftrightarrow$\textup{(\ref{cond:generation})}.
\end{proof}

\section{Applications of Theorem~\ref{thm:regular-comparison}}
\label{sec:applications}

In Subsection~\ref{sec:generation}, we construct the algebra $B$
and use its injective resolution to prove
Corollary~\ref{cor:pure-injective-generation}.
In Subsection~\ref{sec:realization}, we use the same algebra
and resolution to prove Corollary~\ref{cor:separation-example}.

\subsection{Generation by pure injective objects}
\label{sec:generation}

We first record two Ext calculations needed for the construction below.

\begin{lem}\label{lem:tensor-ext}
Let $k$ be a field, let $A,C$ be $k$-algebras, and set
$S=A\otimes_k C$. Let $M$ be an $A$-module, $L$ a
$C$-module, and $N$ an $S$-module. Suppose that
$\Ext_A^q(M,N)=0\;(0\leq q<n)$
for some integer $n\geq0$. Then
$\Ext_S^m(M\otimes_k L,N)=0\;(0\leq m<n)$
and
$\Ext_S^n(M\otimes_k L,N) \cong\Hom_C\bigl(L,\Ext_A^n(M,N)\bigr)$.
The $C$-action on $\Ext_A^n(M,N)$ is induced by its action on
$N$.
\end{lem}
\begin{proof}
Choose an injective $S$-resolution $N\to I^\bullet$ and
write $W^\bullet=\Hom_A(M,I^\bullet)$. Since $S$ is free as a
right $A$-module, restriction preserves injectives, so
$H^q(W)=\Ext_A^q(M,N)$. Each $W^j$ is injective over $C$:
the adjunction
$\Hom_C\bigl(L',\Hom_A(M,I^j)\bigr) \cong\Hom_S(M\otimes_k L',I^j)$
and exactness of $M\otimes_k-$ prove this assertion.

Let $K^j=\im d_W^{j-1}$, with $K^0=0$. Since
$H^j(W)=0$ for $j<n$, the sequences
$$
 0\longrightarrow K^j\longrightarrow W^j
   \longrightarrow K^{j+1}\longrightarrow0
 \;(0\leq j<n)
$$
are exact. They split successively: starting with $K^0=0$, injectivity of $K^j$ splits its
inclusion, and $K^{j+1}$ is an injective direct summand of $W^j$.
Thus all $K^j$ with $j\leq n$ are injective. Applying
$\Hom_C(L,-)$ preserves exactness in degrees less than $n$.
Moreover, $K^n\subseteq Z^n(W)$ splits, and
$W^{n-1}\to K^n$ splits when $n>0$. Therefore
$H^n\Hom_C(L,W^\bullet) \cong\Hom_C(L,H^n(W))$.
For $n=0$, the same isomorphism follows because Hom preserves kernels.
Finally,
$\Hom_C(L,W^\bullet)\cong\Hom_S(M\otimes_k L,I^\bullet)$,
which proves both assertions.
\end{proof}

\begin{lem}\label{lem:extension-orthogonal}
Let $B\to T$ be a ring homomorphism and let $N$ be a
$B$-module. If
$\Ext_B^r(T,N)=0\;(r\geq0)$,
then $\Ext_B^r(M,N)=0$ for every $T$-module $M$ and every
$r\geq0$, where $M$ is restricted to $B$.
\end{lem}
\begin{proof}
Choose an exact sequence of $T$-modules
$0\longrightarrow K\longrightarrow F\longrightarrow M \longrightarrow0$
with $F\cong T^{(I)}$ free. For every $r\geq0$,
$\Ext_B^r(F,N)\cong\prod_I\Ext_B^r(T,N)=0$.
The inclusion
$\Hom_B(M,N)\hookrightarrow\Hom_B(F,N)$ proves the assertion in
degree zero for all $T$-modules. For $r>0$, argue by induction
on $r$, simultaneously for all $T$-modules. The exact sequence
$$
 \Ext_B^{r-1}(K,N)\longrightarrow\Ext_B^r(M,N)
   \longrightarrow\Ext_B^r(F,N)
$$
has zero end terms, so its middle term is zero as well.
\end{proof}

\begin{prop}\label{prop:ext-vanishing}
For every field $k$ there is a faithfully flat inclusion $B\to T$
of commutative von Neumann regular $k$-algebras, both of
$k$-dimension $\aleph_\omega$, such that
$\Ext_B^m(T,B)=0$ for every $m\geq0$.
Moreover, $\Ext_B^m(E,B)=0$ for every injective $B$-module $E$
and every $m\geq0$.
\end{prop}

\begin{proof}
We construct $B\subseteq T=\bigcup_{n<\omega}T_n$ so that
$\Ext_B^m(T_n,B)=0$ for $0\leq m<n$.
The telescope exact sequence then proves
$\Ext_B^m(T,B)=0$ for every $m\geq0$.
Set $\kappa_n=\aleph_n$. For each $n\geq0$, choose a set
$I_n=I_n^+\sqcup I_n^-$ with both parts of cardinality
$\kappa_n$. Let $C_n$ be the $k$-algebra of functions on
$I_n$ which are constant off a finite set. In $k^{I_n}$, let
$q_n$ be the characteristic function of $I_n^+$ and set
$C'_n=C_n[q_n]$. For $i\in I_n$, write $e_i\in C_n$
for the characteristic function of $\{i\}$. Thus
$e_iq_n=e_i$ if $i\in I_n^+$, and $e_iq_n=0$ if $i\in I_n^-$.
Choose a set $K_n$ of cardinality $\kappa_n$, and set
$$
  B_n=\bigotimes_{t\in K_n} C_n,\;
  D_n=\bigotimes_{t\in K_n} C'_n,\;
  B=\bigotimes_{n<\omega}B_n,\; T=\bigotimes_{n<\omega}D_n .
$$
All tensor products in this construction are over $k$.
An infinite tensor product means the direct limit of the
finite tensor products along the unit maps. The maps $B_n\to
D_n$ and $B\to T$ are monic: this holds for the finite
tensor products over a field and passes to the direct limits.
Each ring is generated by idempotents over $k$, hence is a
directed union of finite products of $k$ and is von Neumann
regular. Since every $B$-module is flat, the quotient $T/B$ is flat,
so $B\to T$ is pure. The $B$-module $T$ is flat as well;
purity makes $M\to T\otimes_B M$ monic for every $B$-module
$M$, so the extension is faithfully flat.
A basis of each factor containing its unit yields a basis of
each tensor product consisting of tensors with finitely many
factors different from $1$. Since
$\dim_k C_n=\dim_k C'_n=\kappa_n$, we have
$\dim_k B_n=\dim_k D_n=\kappa_n$ and
$\dim_k B=\dim_k T=\aleph_\omega$.

Write $A_n=\bigotimes_{j<n}B_j$ and
$F_n=\bigotimes_{j<n}D_j$, where $A_0=F_0=k$. Set
$T_n=B\otimes_{A_n} F_n =F_n\otimes_k\bigotimes_{j\geq n} B_j $.
These subalgebras satisfy $T_n\subseteq T_{n+1}$ and
$T=\bigcup_{n<\omega}T_n$. Write $N=B$.
For $n>0$, the inequality $\kappa_{n-1}<\kappa_n$ ensures
that a given cocycle on the free resolution chosen below
uses fewer than $\kappa_n$ tensor factors at level $n$.
We can therefore choose an index $t\in K_n$ whose factor
is not used by the cocycle, and compare evaluations at
points in $I_n^+$ and $I_n^-$.

We first show that, for every $q\geq0$,
$\Hom_{B_n}\bigl(D_n,\Ext_{A_n}^q(F_n,N)\bigr)=0$.
For $n>0$, $A_n$ and $F_n$ have $k$-dimension at most
$\kappa_{n-1}$. Choose a free resolution $P_\bullet\to F_n$
over $A_n$, with each free basis of cardinality at most
$\kappa_{n-1}$. This is possible by taking $k$-bases of
the successive kernels as module generators. For $n=0$ use
the resolution $P_0=k$, $P_q=0$ for $q>0$. In either case
the basis of $P_q$ has cardinality less than $\kappa_n$,
and the union of that many finite sets has cardinality less
than $\kappa_n$.

Let $V=\Ext_{A_n}^q(F_n,N)$, and suppose that
$f:D_n\to V$ is a nonzero $B_n$-linear map. Choose $d\in D_n$
with $f(d)\ne0$. Replacing $f$ by the map $z\mapsto f(dz)$,
we may assume that $u=f(1)\ne0$.
Choose a cocycle for $u$ in $\Hom_{A_n}(P_q,N)$.
Each value on a free basis element is a finite sum of pure
tensors in $B$, with only finitely many factors different from $1$.
The union of the factors from
$K_n$ used by these values has cardinality less than
$\kappa_n$. Choose $t\in K_n$ outside that union.
Denote by $B_{\widehat t}$ the tensor product defining $B$ with
this one copy of $C_n$ omitted. For $i\in I_n$, evaluating this
factor at $i$ gives an $A_n$-algebra retraction
$\pi_i:B\to B_{\widehat t}$. Let
$\pi_{i,*}:V\longrightarrow\Ext_{A_n}^q(F_n,B_{\widehat t})$
be the map induced by $\pi_i$. Since the cocycle for $u$ does not use
the $t$-th factor, all $\pi_{i,*}(u)$ are the same class
$u_0$. The map on Ext induced by $B_{\widehat t}\hookrightarrow B$
sends $u_0$ to $u$, so $u_0\ne0$.

Let $q_t$ denote the copy of $q_n$ in the $t$-th factor
of $D_n$, and set $v=f(q_t)$. Choose a cocycle for $v$.
Express its values on the free basis as finite sums of pure tensors. The functions occurring in the $t$-th factor
are constant outside finite subsets of $I_n$. The union
of these exceptional subsets again has cardinality less
than $\kappa_n$. Choose $i\in I_n^+$ and $j\in I_n^-$
outside that union. The chosen cochain for $v$ has the same image under
evaluation at $i$ and $j$, so
$\pi_{i,*}(v)=\pi_{j,*}(v)$.
On the other hand, $e_i v=e_i u$ and $e_j v=0$, where the
idempotents are in the $t$-th factor of $B_n$. Evaluation
at $i$ sends $e_i$ to $1$, and evaluation at $j$ sends
$e_j$ to $1$. Thus $\pi_{i,*}(v)=u_0$ and
$\pi_{j,*}(v)=0$, a contradiction.
This proves the asserted Hom vanishing.

We claim that $\Ext_B^m(T_n,N)=0$ for $0\leq m<n$.
Induction from $A_n$ to $B$ is exact and preserves
free modules, so
$\Ext_B^m(T_n,N)\cong \Ext_{A_n}^m(F_n,N)$.
The case $n=0$ is vacuous. Suppose the claim holds for $n$.
For the induction step, Lemma~\ref{lem:tensor-ext} applies with
$A=A_n,\; C=B_n,\; M=F_n,\; L=D_n$.
Since $A_{n+1}=A_n\otimes_k B_n$ and
$F_{n+1}=F_n\otimes_k D_n$, the required Ext groups vanish in
degrees less than $n$. In degree $n$, we obtain
$$
 \Ext_{A_{n+1}}^n(F_{n+1},N)
 \cong\Hom_{B_n}\bigl(D_n,\Ext_{A_n}^n(F_n,N)\bigr)=0
$$
by the Hom vanishing proved above. Induction proves the claim.

We now pass to $T$. For any fixed $m\geq0$, discard the first $m+1$ terms
of the sequence $T_n$. The remaining sequence has the
same direct limit. Its telescope exact sequence is
$$
  0\longrightarrow\bigoplus_{n\geq m+1}T_n
   \xrightarrow{\,1-\mathrm{shift}\,}\bigoplus_{n\geq m+1}T_n
   \longrightarrow T\longrightarrow0.
$$
Here the shift sends each summand $T_n$ to $T_{n+1}$ by
the inclusion. The Ext groups into $N$ from either sum vanish in degrees
$m$ and $m-1$, with degree $-1$ interpreted as zero.
Here Ext takes sums in its first variable to products.
The long exact Ext sequence therefore implies
$\Ext_B^m(T,N)=0$.

By Lemma~\ref{lem:extension-orthogonal},
$\Ext_B^m(M,N)=0$ for every $T$-module $M$ and every $m\geq0$.
If $E$ is injective over $B$, purity makes
$E\to T\otimes_B E$ monic, and this monomorphism splits over
$B$. Since $T\otimes_B E$ is a $T$-module,
$\Ext_B^m(E,B)=0$ for every $m\geq0$.
\end{proof}

\begin{proof}[\bf Proof of Corollary~\ref{cor:pure-injective-generation}]
Take $B$ from Proposition~\ref{prop:ext-vanishing} and an
injective resolution $U$ of $B$, placed in nonnegative degrees.
Then $H^0(U)=B\ne0$ and $H^j(U)=0$ for $j\ne0$.
For every injective $B$-module $E$, the vanishing of
$\Ext_B^m(E,B)$ for all $m\ge0$ implies that
$\Hom_B(E,U)$ is exact. Thus $U$ violates
condition~\textup{(\ref{cond:exactness})} of
Theorem~\ref{thm:regular-comparison}.
Since $B$ is von Neumann regular, it is Ding--Chen.
The theorem therefore implies
$\Loc(\Inj(B))\ne\mathbf D(B)$.
Generation by injective modules is equivalent to generation
by pure injective objects in $\mathbf D(B)$, by
\cite[Remark 6.12]{BHM26}, since $B$ is von Neumann regular.
Hence $\Loc(\mathcal P)\ne\mathbf D(B)$.
\end{proof}

\begin{rem}\label{rem:faithfully-flat-extension}
Let $B\to T$ be the faithfully flat extension in
Proposition~\ref{prop:ext-vanishing}. For every $n\in\mathbb Z$,
we have $\Hom_{\mathbf D(B)}(T,B[n])=0$.
The objects $X\in\mathbf D(B)$ such that
$\Hom_{\mathbf D(B)}(X,B[n])=0$ for every $n\in\mathbb Z$
form a localizing subcategory containing $T$ but not $B$.
Thus $B\notin\Loc(T)$.
Corollary~\ref{cor:pure-injective-generation} also follows
from this fact and \cite[Theorem 6.10]{BHM26}, since
$B\to T$ is faithfully flat.

The extension $B\to T$ is also not descendable. Otherwise, $B$ would
belong to the thick tensor ideal generated by $T$
(see \cite[Definition 3.18]{Mat16}), hence to $\Loc(T)$,
since the latter is a tensor ideal, contrary to the preceding
observation. Aoki constructs a nondescendable faithfully flat
map between Boolean rings (see \cite[Theorem D]{Aok26}).
\end{rem}

\subsection{A coherent ring over which Gorenstein and Ding injectivity differ}
\label{sec:realization}

\begin{samepage}
To prove the local assertion in Corollary~\ref{cor:separation-example}, we
choose an Ext class whose localization is nonzero at every prime ideal of $R$.

\begin{prop}\label{prop:faithful-class}
Let $B$ be a commutative von Neumann regular ring, and write
$R=B[x,y]/(x,y)^2$ and $\mathfrak r=(x,y)$. Suppose that
$U$ is a complex of injective $B$-modules for which
$\Hom_B(E,U)$ is exact for every injective
$B$-module $E$, and that $[w]\in H^0(U)$ satisfies
$\ann_B[w]=0$. For the periodic
module $G$ of Proposition~\ref{prop:periodic-transfer},
there is a class $\xi\in\Ext_R^1(F_B,G)$ such that
$$
 \ann_R\xi=\mathfrak r,\;
 B=R/\mathfrak r\lhook\joinrel\longrightarrow
       \Ext_R^1(F_B,G),\; b\longmapsto b\xi.
$$
Moreover, $\Did_R G=\infty$ and
$\Gid_{R_{\mathfrak p}}G_{\mathfrak p}=\infty$ for every
$\mathfrak p\in\Spec R$.
\end{prop}
\end{samepage}
\begin{proof}
Let $\xi$ be the image of $[w]$ under the monomorphism
of Proposition~\ref{prop:periodic-transfer}. This map is $B$-linear and
injective, and its image is annihilated by $\mathfrak r$.
Therefore
$\ann_R\xi=\mathfrak r+\ann_B[w]=\mathfrak r$.
For every $\mathfrak p\in\Spec R$,
$\ann_R\xi=\mathfrak r\subseteq\mathfrak p$, so
$\xi_{\mathfrak p}\ne0$.
The periodic injective sequence induces $R$-linear
dimension-shifting isomorphisms
$\Ext_R^n(F_B,G)\cong\Ext_R^1(F_B,G)\;(n\ge1)$.
The ring $R$ is coherent and $F_B$ is finitely presented
by Lemma~\ref{lem:finite-detector}. Thus $F_B$ has a
resolution by finitely generated free modules. Localizing
that resolution, we obtain
$$
 \Ext_{R_{\mathfrak p}}^n((F_B)_{\mathfrak p},G_{\mathfrak p})
 \cong\Ext_R^n(F_B,G)_{\mathfrak p}\ne0
 \;(n\ge1).
$$
The module $(F_B)_{\mathfrak p}$ is injective over $R_{\mathfrak p}$
by the same lemma. A finite Gorenstein injective coresolution
of $G_{\mathfrak p}$ over $R_{\mathfrak p}$ would make these Ext groups zero in all
sufficiently large degrees, a contradiction. Similarly, a
finite Ding injective coresolution of $G$ over $R$ would make
$\Ext_R^n(F_B,G)$ zero for all sufficiently large $n$,
since $F_B$ is absolutely pure. Consequently,
$\Did_R G=\infty$ and
$\Gid_{R_{\mathfrak p}}G_{\mathfrak p}=\infty$ for every
$\mathfrak p\in\Spec R$.
\end{proof}

\begin{samepage}
\begin{proof}[\bf Proof of Corollary~\ref{cor:separation-example}]
Take the algebra $B$ and its injective resolution $U$ used
in the proof of Corollary~\ref{cor:pure-injective-generation}.
Set $R=B\otimes_k\Lambda_k$.
Then $\dim_kR=\aleph_\omega$.
By Lemma~\ref{lem:finite-detector}, $R$ is coherent
and $F_B=\Hom_B(R,B)$ is finitely presented and
absolutely pure. By Proposition~\ref{prop:single-source},
$\DI(R)=\Inj(R)$.

Since $B$ is von Neumann regular, it is Ding--Chen.
The complex $U$ violates condition~\textup{(\ref{cond:exactness})}
of Theorem~\ref{thm:regular-comparison}, so $\DI(R)\subsetneq\GI(R)$.
Let $G$ be the module constructed from $U$ in
Proposition~\ref{prop:periodic-transfer}. Then $G$ is strongly
Gorenstein injective and is not Ding injective.
Under $H^0(U)\cong B$, the class corresponding to $1$ has
zero annihilator in $B$. It follows from
Proposition~\ref{prop:faithful-class} that $\Did_RG=\infty$ and
$\Gid_{R_{\mathfrak p}}G_{\mathfrak p}=\infty$
for every prime ideal $\mathfrak p$ of $R$.
Since $H^0(U)\cong B$ and $H^j(U)=0$ for $j\ne0$,
the module $W$ in Proposition~\ref{prop:periodic-transfer}
is isomorphic to $B^{\mathbb N}$. Thus
$\Ext_R^n(F_B,G)\cong(B^{\mathbb N})^3\;(n\ge1)$,
where the $R$-action on the right factors through $R\to B$.

Every finite set of coefficients in $B$ lies in a
$k$-subalgebra which is a finite product of copies of
$k$. Tensoring these subalgebras with $\Lambda_k$
expresses $R$ as a directed union of finite products
of $\Lambda_k$.
Every residue field of $B$ is $k$, since the generating
idempotents take values zero or one in a field.
Lemma~\ref{lem:finite-detector} therefore identifies
every $R_{\mathfrak p}$ with $\Lambda_k$.
\end{proof}
\end{samepage}

\begin{rem}\label{rem:further-consequences}
Let $R=B\otimes_k\Lambda_k$ and $G$ be as in the proof of
Corollary~\ref{cor:separation-example}.
The ring $R$ is a directed union of finite products of $\Lambda_k$.
These finite stages are Noetherian, so Gorenstein and Ding
injective modules coincide over each of them.
Thus the equality $\GI=\DI$ need not be preserved under
directed unions of rings.

Restriction to $\Lambda_k$ need not preserve Gorenstein
injectivity, although $R$ is free over $\Lambda_k$.
Indeed, $E_0=\Hom_k(\Lambda_k,k)$ is injective over $\Lambda_k$
and $R\otimes_{\Lambda_k}E_0\cong F_B$. Ext adjunction gives
$\Ext_{\Lambda_k}^1(E_0,G_{\Lambda_k})
\cong\Ext_R^1(F_B,G)\ne0$.

The same Ext calculation shows that ${}^{\perp}\GI(R)$ is
not closed under pure submodules. The embedding of $F_B$
into its injective envelope $I$ is pure, since $F_B$ is
absolutely pure. One has $I\in{}^{\perp}\GI(R)$ and
$F_B\notin{}^{\perp}\GI(R)$.
\end{rem}

\begin{cor}\label{cor:strict-three-classes}
There is a commutative coherent ring $R'$ of cardinality
$\aleph_\omega$ with
$\Inj(R')\subsetneq\DI(R')\subsetneq\GI(R')$.
\end{cor}
\begin{proof}
We first show that
$\Inj(S)\subsetneq\DI(S)$ for every nonzero commutative ring
$A$, where $S=(A[x,y]/(x,y)^2)[\varepsilon]/(\varepsilon^2)$.
Set $L=S/\varepsilon S$.
Since $\ann_{S}(\varepsilon)=\varepsilon S$, the
module $L$ has a free resolution
$$
 \cdots\xrightarrow{\varepsilon}S
 \xrightarrow{\varepsilon}S\longrightarrow L
 \longrightarrow0.
$$
Consequently, for every $S$-module $Q$,
$\Ext_{S}^1(L,Q) \cong\ker(\varepsilon:Q\to Q)/\varepsilon Q$.
Choose a nonzero injective $S$-module $I$.
The one-periodic complex $J$ with term $I$ and
differential $\varepsilon$ is exact, since
$\Ext_{S}^1(L,I)=0$.
If $Q$ is absolutely pure, then $\Ext_{S}^1(L,Q)$
vanishes because $L$ is finitely presented.
Thus the one-periodic complex with term $Q$ and
differential $\varepsilon$ is exact.
Applying the exact functor $\Hom_{S}(-,I)$ proves that
$\Hom_{S}(Q,J)$ is exact:
precomposition with $\varepsilon$ on $Q$ agrees with
postcomposition with $\varepsilon$ on $I$.
Hence $N=\varepsilon I=\ker(\varepsilon:I\to I)$ is
Ding injective.
It is nonzero, since otherwise exactness of $J$ would
force $I=0$.
As $\varepsilon N=0$, the same resolution gives
$\Ext_{S}^1(L,N)\cong N\ne0$.
Thus $N$ is not injective.

Take a finite field $k$, and let $R=B[x,y]/(x,y)^2$
and $G$ be as in the proof of Corollary~\ref{cor:separation-example}.
Write $F_B=\Hom_B(R,B)$.
Set
$$
 R'=R[\varepsilon]/(\varepsilon^2)
   \cong \bigl(B[\varepsilon]/(\varepsilon^2)\bigr)[x,y]/(x,y)^2.
$$
Choose a complete injective resolution $J$ of $G$.
The functor $\Hom_R(R',-)$ preserves injectives by adjunction
and is exact because $R'$ is free over $R$.
For every injective $R'$-module $Q$, freeness also makes
$Q_R$ injective over $R$. By adjunction,
$\Hom_{R'}(Q,\Hom_R(R',J))\cong\Hom_R(Q_R,J)$ is exact.
Thus $\Hom_R(R',J)$ is totally acyclic, and
$G'=\Hom_R(R',G)$ is Gorenstein injective over $R'$.

Let $F'=\Hom_R(R',F_B)$.
The $R'$-module $F'$ is absolutely pure: restriction takes
finitely presented $R'$-modules to finitely presented
$R$-modules, so this follows from Ext adjunction and
absolute purity of $F_B$ over $R$.
Since $(F')_R\cong F_B^2$, the same adjunction gives
$$
 \Ext_{R'}^1(F',G')
 \cong\Ext_R^1((F')_R,G)
 \cong\Ext_R^1(F_B,G)^2\ne0.
$$
Thus $G'$ is not Ding injective.
Taking $A=B$ in the preceding construction proves
$\Inj(R')\subsetneq\DI(R')$.
Finally, $R'$ is finite free over $R$, so it is
coherent by the kernel argument of
Lemma~\ref{lem:finite-detector}, and
$\lvert R'\rvert=\aleph_\omega$.
\end{proof}

\smallskip
{\bf Acknowledgements.} The first author was supported by the National Natural Science
Foundation of China (Grant No.~11901412).
The second author was supported by the National Natural Science
Foundation of China (Grant No.~12571035).

\smallskip
{\bf The Use of AI.} Generative artificial intelligence tools were used to assist with
literature searches, the formulation and checking of mathematical
arguments, language and structural revision, and LaTeX preparation.
The authors are responsible for the mathematical content,
references, and final text.

\vspace{3mm}
\noindent\textbf{Guocheng Dai}\\
School of Mathematical Sciences, Sichuan Normal University, Chengdu 610066, P. R. China.\\
Email: \textsf{ggcdai@163.com}\\[1mm]
\textbf{Jiangsheng Hu}\\
School of Mathematics, Hangzhou Normal University, Hangzhou 311121, P. R. China.\\
Email: \textsf{hujs@hznu.edu.cn}\\[1mm]

\bigskip
\begin{thebibliography}{99}
\bibitem{Aok26} K. Aoki,
\emph{On cohomology of locally profinite sets},
Proc. Lond. Math. Soc. (3) \textbf{132} (2026), no.~5, e70159.
\bibitem{BCE} S. Bazzoni, M. Cort\'es-Izurdiaga and S. Estrada,
\emph{Periodic modules and acyclic complexes},
Algebr. Represent. Theory \textbf{23} (2020), no.~5, 1861--1883.
\bibitem{BM07} D. Bennis and N. Mahdou,
\emph{Strongly Gorenstein projective, injective, and flat modules},
J. Pure Appl. Algebra \textbf{210} (2007), no.~2, 437--445.
\bibitem{BHM26} S. Breaz, M. Hrbek and G. C. Modoi,
\emph{Silting, cosilting and extensions of commutative rings},
Math. Proc. Cambridge Philos. Soc. \textbf{180} (2026), no.~2, 363--388.
\bibitem{Cou06} F. Couchot,
\emph{Localization of injective modules over valuation rings},
Proc. Amer. Math. Soc. \textbf{134} (2006), no.~4, 1013--1017.
\bibitem{EJ95} E. E. Enochs and O. M. G. Jenda,
\emph{Gorenstein injective and projective modules},
Math. Z. \textbf{220} (1995), 611--633.
\bibitem{Gil10} J. Gillespie,
\emph{Model structures on modules over Ding--Chen rings},
Homology Homotopy Appl. \textbf{12} (2010), no.~1, 61--73.
\bibitem{Gil17} J. Gillespie,
\emph{On Ding injective, Ding projective, and Ding flat modules and complexes},
Rocky Mountain J. Math. \textbf{47} (2017), no.~8, 2641--2673.
\bibitem{Mao23} L. Mao,
\emph{Ding projective and Ding injective modules over
trivial ring extensions},
Czechoslovak Math. J. \textbf{73} (2023), no.~3, 903--919.
\bibitem{MD08} L. Mao and N. Ding,
\emph{Gorenstein FP-injective and Gorenstein flat modules},
J. Algebra Appl. \textbf{7} (2008), no.~4, 491--506.
\bibitem{Mat16} A. Mathew,
\emph{The Galois group of a stable homotopy theory},
Adv. Math. \textbf{291} (2016), 403--541.
\bibitem{Pos17} L. Positselski,
\emph{Coherent rings, fp-injective modules, dualizing complexes,
and covariant Serre--Grothendieck duality},
Selecta Math. (N.S.) \textbf{23} (2017), no.~2, 1279--1307.
\bibitem{Ric19} J. Rickard,
\emph{Unbounded derived categories and the finitistic dimension conjecture},
Adv. Math. \textbf{354} (2019), 106735.
\bibitem{Spa88} N. Spaltenstein,
\emph{Resolutions of unbounded complexes},
Compositio Math. \textbf{65} (1988), no.~2, 121--154.
\bibitem{WLY21} J. Wang, Z. Liu and G. Yang,
\emph{Gillespie's questions and Grothendieck duality},
C. R. Math. \textbf{359} (2021), no.~5, 593--607.
\bibitem{ZhaPeriodic} C. Zhang,
\emph{Commutative coherent rings with
$\mathcal{GP}\not\subseteq\mathcal{GF}$},
preprint, arXiv:2608.17748v3 (2026).
\end{thebibliography}
\end{document}